\documentclass{article}

\usepackage[dvipsnames]{xcolor}

\usepackage[marginparsep=1cm,top=4cm,bottom=3cm,inner=2.5cm,outer=5.5cm,headsep=30pt,a4paper]{geometry} 

\usepackage{mparhack}

\usepackage[parfill]{parskip}    

\usepackage[export]{adjustbox}  

\usepackage{stmaryrd}

\usepackage{datetime}

\newcommand{\BM}[1]{{\mbox{\boldmath$#1$}}}

\newcommand{\BMR}[1]{{\mbox{\boldmath$\mathrm{#1}$}}}

\usepackage[utf8]{inputenc} 
\usepackage[T1]{fontenc} 
\usepackage{gfsartemisia-euler}

\usepackage{pifont}

\usepackage{csquotes}

\usepackage{fontawesome}
\usepackage{bbding}

\usepackage{enumitem} 
\setlist{nolistsep} 

\usepackage{booktabs} 

\usepackage{manfnt}

\usepackage[dvipsnames]{xcolor} 

\usepackage{fontawesome}

\definecolor{ocretriad}{RGB}{243,90,0} 
\definecolor{darkbluetriad}{RGB}{27,18,178} 
\definecolor{purpletriad}{RGB}{14,0,255} 
\definecolor{darkpurpletriad}{RGB}{10,0,178} 
\definecolor{greentriad}{RGB}{50,204,20} 
\definecolor{darkgreentriad}{RGB}{36,178,9} 

\definecolor{yellowbckgr}{RGB}{255,244,168} 
\definecolor{yellowline}{RGB}{244,221,41} 
\definecolor{yellowtext}{RGB}{200,165,0} 

\definecolor{ocre}{RGB}{243,90,0} 
\definecolor{goldline}{RGB}{194, 152, 18}
\definecolor{goldtext}{RGB}{157, 123, 15}
\definecolor{sunny}{RGB}{252, 180, 80}

\definecolor{greentext}{RGB}{30, 122, 12}
\definecolor{greenline}{RGB}{25, 107, 49}
\definecolor{greenbckgr}{RGB}{111, 220, 80}

\definecolor{bluetext}{RGB}{47, 96, 121}
\definecolor{bluebckgr}{RGB}{57, 188, 214}
\definecolor{blueline}{RGB}{56, 115, 145}

\definecolor{purpleline}{RGB}{91,55,156} 
\definecolor{purpletext}{RGB}{64,6,127} 
\definecolor{purplebckgr}{RGB}{140,90,255} 

\definecolor{redline}{RGB}{199,73,51} 
\definecolor{nicered}{RGB}{209, 48, 0}
\definecolor{redbckgr}{RGB}{209, 48, 0}

\usepackage{graphicx} 
\graphicspath{{../Pictures/}} 

\usepackage{tikz} 
\usetikzlibrary{shapes.arrows}

\usepackage{fancyhdr} 

\fancypagestyle{plain}{\fancyhead{}} 

\makeatletter
\renewcommand{\cleardoublepage}{
\clearpage\ifodd\c@page\else
\hbox{}
\vspace*{\fill}
\thispagestyle{empty}
\newpage
\fi}

\usepackage{amsmath,amsfonts,amssymb,amsthm} 
\usepackage{mathtools}
\usepackage[mathscr]{eucal}

\newtheoremstyle{ocrenumbox}
{1.4em}
{}
{\normalfont}
{}
{\sffamily\bfseries\color{ocre}}
{\;}
{.5em}
{\sffamily\bfseries\color{ocre}{{\large\thmnumber{\@ifnotempty{#1}{}}\@upn{#2}}\nobreakspace{\large\thmname{#1}}}
\nobreakspace\nobreakspace\nobreakspace
\thmnote{\the\thm@notefont\sffamily\bfseries\itshape\color{Black}{\faHandORight\;\;#3\;\;\;}}
\smallskip\smallskip\newline}

\newtheoremstyle{greennumbox}
{1.4em}
{}
{\normalfont}
{}
{\bf\sffamily\color{greentext}}
{\;}
{.5em}
{\sffamily\bfseries\color{greentext}{{\large\thmnumber{\@ifnotempty{#1}{}}\@upn{#2}}\nobreakspace\thmname{#1}}
\nobreakspace\nobreakspace\nobreakspace
\thmnote{\the\thm@notefont\sffamily\bfseries\itshape\color{Black}{\faHandORight\;\;#3\;\;\;}}
\smallskip\smallskip\newline}

\newtheoremstyle{blacknumbox} 
{1.4em}
{}
{\normalfont}
{}
{\bf\sffamily}
{\;}
{.5em}
{\sffamily{{\large\thmnumber{\@ifnotempty{#1}{}}\@upn{#2}}\nobreakspace\thmname{#1}}
\nobreakspace\nobreakspace\nobreakspace
\thmnote{\the\thm@notefont\sffamily\bfseries\itshape\color{Black}{\faHandORight\;\;#3\;\;\;}}
\smallskip\smallskip\newline}

\newtheoremstyle{yellownumbox}
{1.4em}
{}
{\normalfont}
{}
{\bf\sffamily\color{yellowtext}}
{\;}
{.5em}
{\sffamily\bfseries\color{yellowtext}{{\large\thmnumber{\@ifnotempty{#1}{}}\@upn{#2}}\nobreakspace\thmname{#1}}
\nobreakspace\nobreakspace\nobreakspace
\thmnote{\the\thm@notefont\sffamily\bfseries\itshape\color{Black}{\faHandORight\;\;#3\;\;\;}}
\smallskip\smallskip\newline} 

\newtheoremstyle{bluenumbox}
{1.4em}
{}
{\normalfont}
{}
{\bf\sffamily\color{bluetext}}
{\;}
{.5em}
{\sffamily\bfseries\color{bluetext}{{\large\thmnumber{\@ifnotempty{#1}{}}\@upn{#2}}\nobreakspace\thmname{#1}}
\nobreakspace\nobreakspace\nobreakspace
\thmnote{\the\thm@notefont\sffamily\bfseries\itshape\color{Black}{\faHandORight\;\;#3\;\;\;}}
\smallskip\smallskip\newline}

\newtheoremstyle{rednumbox}
{1.4em}
{}
{\normalfont}
{}
{\bf\color{nicered}}
{\;}
{.5em}
{\sffamily\bfseries\color{nicered}{{\large\thmnumber{\@ifnotempty{#1}{}}\@upn{#2}}\nobreakspace\thmname{#1}}
\nobreakspace\nobreakspace\nobreakspace
\thmnote{\the\thm@notefont\sffamily\bfseries\itshape\color{Black}{\faHandORight\;\;#3\;\;\;}}
\smallskip\smallskip\newline}

\newtheoremstyle{purplenumbox}
{1.4em}
{}
{\normalfont}
{}
{\bf\sffamily\color{purpletext}}
{\;}
{.5em}
{\sffamily\bfseries\color{purpletext}{{\large\thmnumber{\@ifnotempty{#1}{}}\@upn{#2}}\nobreakspace\thmname{#1}}
\nobreakspace\nobreakspace\nobreakspace
\thmnote{\the\thm@notefont\sffamily\bfseries\itshape\color{Black}{\faHandORight\;\;#3\;\;\;}}
\smallskip\smallskip\newline}

\newtheoremstyle{goldnumbox}
{1.4em}
{}
{\normalfont}
{}
{\bf\sffamily\color{goldtext}}
{\;}
{.5em}
{\sffamily\bfseries\color{goldtext}{{\large\thmnumber{\@ifnotempty{#1}{}}\@upn{#2}}\nobreakspace\thmname{#1}}
\nobreakspace\nobreakspace\nobreakspace
\thmnote{\the\thm@notefont\sffamily\bfseries\itshape\color{Black}{\faHandORight\;\;#3\;\;\;}}
\smallskip\smallskip\newline}

\newtheoremstyle{blacknonumbox}
{1.4em}
{}
{\normalfont}
{} 
{\bf\sffamily}
{:}
{.5em}
{
\nobreakspace{{\large\thmname{#1}}}
\nobreakspace\nobreakspace\nobreakspace
\thmnote{\the\thm@notefont\sffamily\bfseries\itshape\color{Black}{\faHandORight\;\;#3\;\;\;}}
\smallskip\smallskip\newline}

\newtheoremstyle{ocrenum}
{1.5em}
{}
{\normalfont}
{}
{\bf\sffamily\color{ocre}}
{\;}
{.5em}
{\sffamily\bfseries\color{ocre}{{\large\thmnumber{\@ifnotempty{#1}{}}\@upn{#2}}\nobreakspace\thmname{#1}}
\nobreakspace\nobreakspace\nobreakspace
\thmnote{\the\thm@notefont\sffamily\bfseries\itshape\color{Black}{\faHandORight\;\;#3\;\;\;}}
\smallskip\smallskip\newline} 
\makeatother

\newcounter{dummy} 
\numberwithin{dummy}{section}

\theoremstyle{ocrenumbox}
\newtheorem{theoremT}[dummy]{Theorem}
\newtheorem{factT}[dummy]{Fact}
\newtheorem{corollaryT}[dummy]{Corollary}
\newtheorem{lemmaT}[dummy]{Lemma}
\newtheorem{AlgT}[dummy]{Algorithm}

\theoremstyle{greennumbox}
\newtheorem{exerciseT}[dummy]{Exercise}
\newtheorem{TYCT}[dummy]{Test Your Comprehension}

\theoremstyle{bluenumbox}
\newtheorem{definitionT}[dummy]{Definition}
\newtheorem{TrmT}[dummy]{Terminology}
\newtheorem{NtnT}[dummy]{Notation}

\theoremstyle{goldnumbox}
\newtheorem{exampleT}[dummy]{Example}
\newtheorem{examplesT}[dummy]{Examples}

\theoremstyle{rednumbox}
\newtheorem{ComT}[dummy]{Comment}
\newtheorem{SynT}[dummy]{Synopsis}

\theoremstyle{purplenumbox}
\newtheorem{ObsT}[dummy]{Observation}

\theoremstyle{blacknonumbox}

\theoremstyle{blacknumbox}
\newtheorem{ProbT}[dummy]{Problem}
\newtheorem{BProbT}[dummy]{Bonus Problem}

\RequirePackage[framemethod=TikZ]{mdframed} 

\newmdenv[skipabove=7pt,
skipbelow=7pt,
backgroundcolor=sunny!25,
linecolor=ocre,
roundcorner=10pt,
innerleftmargin=10pt,
innerrightmargin=10pt,
innertopmargin=0pt,
leftmargin=0cm,
rightmargin=0cm,
innerbottommargin=8pt,
linewidth=2pt,
roundcorner=10pt,
splittopskip=1.5\baselineskip]{tBox}

\mdfdefinestyle{FormulaStyle}{skipabove=23pt,
skipbelow=0pt,
backgroundcolor=sunny!25,
linecolor=ocre,
rightline=true,
leftline=true,
topline=true,
bottomline=true,
linewidth=2pt,
innerleftmargin=10pt,
innerrightmargin=10pt,
innertopmargin=8pt,
innerbottommargin=8pt,
leftmargin=0cm,
rightmargin=0cm,
roundcorner=10pt,
splittopskip=1.5\baselineskip
}

\mdfdefinestyle{GeometryStyle}{skipabove=23pt,
skipbelow=0pt,
backgroundcolor=Gray!10,
linecolor=Black,
rightline=true,
leftline=true,
topline=true,
bottomline=true,
linewidth=2pt,
innerleftmargin=10pt,
innerrightmargin=10pt,
innertopmargin=8pt,
innerbottommargin=8pt,
leftmargin=0cm,
rightmargin=0cm,
roundcorner=10pt,
splittopskip=1.5\baselineskip}

\mdfdefinestyle{ProofSynopsisStyle}{skipabove=23pt,
skipbelow=7pt,
backgroundcolor=black!5,
linecolor=ocre,
rightline=false,
leftline=false,
topline=false,
bottomline=false,
linewidth=2pt,
innerleftmargin=10pt,
innerrightmargin=10pt,
innertopmargin=8pt,
innerbottommargin=8pt,
leftmargin=0cm,
rightmargin=0cm,
splittopskip=1.5\baselineskip
}

\newmdenv[skipabove=7pt,
skipbelow=7pt,
backgroundcolor=OrangeRed!6,
linecolor=redline,
innerleftmargin=10pt,
innerrightmargin=10pt,
innertopmargin=0pt,
leftmargin=0cm,
rightmargin=0cm,
innerbottommargin=8pt,
linewidth=2pt,
roundcorner=10pt,
splittopskip=1.5\baselineskip]{SynBox}

\newmdenv[skipabove=7pt,
skipbelow=7pt,
rightline=false,
leftline=true,
topline=false,
bottomline=true,
linecolor=goldline,
innerleftmargin=10pt,
innerrightmargin=5pt,
innertopmargin=0pt,
innerbottommargin=8pt,
leftmargin=0cm,
rightmargin=0cm,
linewidth=2pt,
roundcorner=10pt,
splittopskip=1.5\baselineskip]{exmBox}	

\newmdenv[skipabove=7pt,
skipbelow=7pt,
rightline=false,
leftline=true,
topline=false,
bottomline=true,
backgroundcolor=greentriad!25,
linecolor=greenline,
innerleftmargin=10pt,
innerrightmargin=5pt,
innertopmargin=0pt,
innerbottommargin=8pt,
leftmargin=0cm,
rightmargin=0cm,
linewidth=2pt,
roundcorner=10pt,
splittopskip=1.5\baselineskip]{eBox}

\newmdenv[skipabove=5pt,
skipbelow=7pt,
rightline=false,
leftline=true,
topline=false,
bottomline=true,
linecolor=black,
innerleftmargin=10pt,
innerrightmargin=5pt,
innertopmargin=0pt,
innerbottommargin=8pt,
leftmargin=0cm,
rightmargin=0cm,
linewidth=2pt,
roundcorner=10pt,
splittopskip=1.5\baselineskip]{ProbBox}

\newmdenv[skipabove=7pt,
skipbelow=7pt,
rightline=false,
leftline=true,
topline=false,
bottomline=true,
backgroundcolor=greentriad!10,
linecolor=greenline,
innerleftmargin=10pt,
innerrightmargin=5pt,
innertopmargin=0pt,
innerbottommargin=10pt,
leftmargin=0cm,
rightmargin=0cm,
linewidth=2pt,
roundcorner=10pt,
splittopskip=1.5\baselineskip]{TYCBox}

\newmdenv[skipabove=7pt,
skipbelow=7pt,
rightline=false,
leftline=true,
topline=false,
bottomline=true,
backgroundcolor=yellowbckgr!50,
linecolor=yellowline,
innerleftmargin=10pt,
innerrightmargin=5pt,
innertopmargin=0pt,
innerbottommargin=10pt,
leftmargin=0cm,
rightmargin=0cm,
linewidth=2pt,
roundcorner=10pt,
splittopskip=1.5\baselineskip]{AlgBox}

\newmdenv[skipabove=7pt,
skipbelow=7pt,
rightline=false,
leftline=true,
topline=false,
bottomline=true,
linecolor=purpleline,
innerleftmargin=10pt,
innerrightmargin=5pt,
innertopmargin=0pt,
innerbottommargin=8pt,
leftmargin=0cm,
rightmargin=0cm,
linewidth=2pt,
roundcorner=10pt,
splittopskip=1.5\baselineskip]{ObsBox}

\newmdenv[skipabove=7pt,
skipbelow=7pt,
rightline=false,
leftline=true,
topline=false,
bottomline=true,
linecolor=blueline,
innerleftmargin=10pt,
innerrightmargin=5pt,
innertopmargin=0pt,
leftmargin=0cm,
rightmargin=0cm,
linewidth=2pt,
innerbottommargin=8pt,
roundcorner=10pt,
splittopskip=1.5\baselineskip]{dBox}

\newmdenv[skipabove=7pt,
skipbelow=7pt,
rightline=false,
leftline=true,
topline=false,
bottomline=true,
linecolor=blueline,
innerleftmargin=10pt,
innerrightmargin=5pt,
innertopmargin=0pt,
leftmargin=0cm,
rightmargin=0cm,
linewidth=2pt,
innerbottommargin=8pt,
roundcorner=10pt,
splittopskip=1.5\baselineskip]{TrmBox}

\newmdenv[skipabove=7pt,
skipbelow=7pt,
rightline=true,
leftline=true,
topline=true,
bottomline=true,
linecolor=ocre,
backgroundcolor=sunny!25,
innerleftmargin=10pt,
innerrightmargin=5pt,
innertopmargin=0pt,
leftmargin=0cm,
rightmargin=0cm,
linewidth=2pt,
innerbottommargin=8pt,
roundcorner=10pt,
splittopskip=1.5\baselineskip]{cBox}

\newmdenv[skipabove=7pt,
skipbelow=7pt,
rightline=false,
leftline=true,
topline=false,
bottomline=true,
linecolor=redline,
innerleftmargin=10pt,
innerrightmargin=5pt,
innertopmargin=0pt,
leftmargin=0cm,
rightmargin=0cm,
linewidth=2pt,
innerbottommargin=8pt,
roundcorner=10pt,
splittopskip=1.5\baselineskip]{ComBox}

\newmdenv[skipabove=7pt,
skipbelow=7pt,
rightline=false,
leftline=true,
topline=false,
bottomline=true,
linecolor=black,
innerleftmargin=10pt,
innerrightmargin=5pt,
innertopmargin=0pt,
innerbottommargin=8pt,
leftmargin=0cm,
rightmargin=0cm,
linewidth=2pt,
roundcorner=10pt,
splittopskip=1.5\baselineskip
]{PfBox}

\newenvironment{Thm}{\begin{tBox}\begin{theoremT}}{\end{theoremT}\end{tBox}}

\newenvironment{Dfn}{\begin{dBox}\begin{definitionT}}{\end{definitionT}\end{dBox}}

\newenvironment{Obs}{\begin{ObsBox}\begin{ObsT}}{\end{ObsT}\end{ObsBox}}
\newenvironment{Cor}{\begin{cBox}\begin{corollaryT}}{\end{corollaryT}\end{cBox}}

\newenvironment{Formula}{\Vm{-1.5}\begin{mdframed}[style=FormulaStyle]}{\end{mdframed}}

\makeatletter

\renewcommand{\@seccntformat}[1]{\llap{\textcolor{black!66}{\csname the#1\endcsname} 
\hspace{1em}}}   

\renewcommand{\section}{\@startsection{section}{1}{\z@}
{-4ex \@plus -1ex \@minus -.4ex}
{1ex \@plus.2ex }
{\normalfont\LARGE\sffamily\bfseries}}

\renewcommand{\subsection}{\@startsection{subsection}{2}{\z@}
{-3ex \@plus -0.1ex \@minus -.4ex}
{0.5ex \@plus.2ex }
{\normalfont\Large\sffamily\bfseries}}

\renewcommand{\subsubsection}{\@startsection{subsubsection}{3}{\z@}
{-2ex \@plus -0.1ex \@minus -.2ex}
{.2ex \@plus.2ex }
{\normalfont\sffamily\bfseries}}                        

\renewcommand\paragraph{\@startsection{paragraph}{4}{\z@}
{-2ex \@plus-.2ex \@minus .2ex}
{.1ex}
{\normalfont\small\sffamily\bfseries}}

\usepackage{fancybox}

\usepackage{sidenotes}

\usepackage{hyperref}
\hypersetup{hidelinks,colorlinks=false,breaklinks=true,urlcolor= ocre,bookmarksopen=false,pdftitle={Title},pdfauthor={Author}}
\usepackage{bookmark}
\bookmarksetup{
open,
numbered,
addtohook={%
\ifnum\bookmarkget{level}=0 
\bookmarksetup{bold}%
\fi
\ifnum\bookmarkget{level}=-1 
\bookmarksetup{color=ocre,bold}%
\fi
}
}

\usepackage{stackengine}
\stackMath

\newcommand{\Fam}[1]{\ensuremath{\mathcal{#1}}}

\newcommand{\Field}[1]{\ensuremath{\mathbb{#1}}}

\newcommand{\C}{\Field{C}}

\newcommand{\R}{\Field{R}}

\newcommand{\F}{\Field{F}}

\newcommand{\Mn}[1]{\ensuremath{\Field{M}_{_{#1}}}}

\newcommand{\Set}[2]{\ensuremath{\left\{\left. \ {#1}\ \right| \ {#2}\ 
\right\} } }

\newcommand{\AND}{\text{\; and \;}}

\newcommand{\Rank}[1]{\ensuremath{\textit{Rank}\left({#1}\right)}}

\newcommand{\Nlty}[1]{\ensuremath{\textit{Nullity}\left({#1}\right)}}

\newcommand{\Span}[1]{\ensuremath{\textit{Span}\left({#1}\right)}}

\newcommand{\Dim}[1]{\ensuremath{\textit{dim}\left({#1}\right)}}
\newcommand{\RREF}[1]{\ensuremath{\textup{RREF}\!\left({#1}\right)}}
\newcommand{\RCEF}[1]{\ensuremath{\textup{RCEF}\!\left({#1}\right)}}



\DeclareFontFamily{U}{matha}{\hyphenchar\font45}
\DeclareFontShape{U}{matha}{m}{n}{
      <5> <6> <7> <8> <9> <10> gen * matha
      <10.95> matha10 <12> <14.4> <17.28> <20.74> <24.88> matha12
      }{}
\DeclareSymbolFont{matha}{U}{matha}{m}{n}
\DeclareFontFamily{U}{mathx}{\hyphenchar\font45}
\DeclareFontShape{U}{mathx}{m}{n}{
      <5> <6> <7> <8> <9> <10>
      <10.95> <12> <14.4> <17.28> <20.74> <24.88>
      mathx10
      }{}
\DeclareSymbolFont{mathx}{U}{mathx}{m}{n}

\DeclareMathSymbol{\obot}         {2}{matha}{"6B}
\DeclareMathSymbol{\bigobot}       {1}{mathx}{"CB}

\makeatletter
\newcommand{\dashline}{\rotatebox[origin=c]{90}{$\dabar@\dabar@\dabar@$}}
\makeatother

\newcommand{\REV}{\ensuremath{\looparrowleft}}

\newcommand{\lra}{\longrightarrow}

\newcommand{\Vm}[1]{\vspace{#1 em}}
\newcommand{\I}{\item}
\newcommand{\LL}{\label}

\newcommand{\BSK}{\bigskip}
\newcommand{\SSK}{\smallskip\smallskip}

\usepackage[]{mdframed}

\newcommand{\RRank}[1]{\ensuremath{\textit{Row Rank}\left({#1}\right)}}
\newcommand{\CRank}[1]{\ensuremath{\textit{Column Rank}\left({#1}\right)}}

\newcommand{\MD}[1]{\begin{Formula}#1\end{Formula}}

\newcommand{\Red}{{\color{red!60!black}red\,}}

\newcommand{\Redim}{{\color{red!60!black}Red}-imension }
\newcommand{\redim}{{\color{red!60!black}red}-imension }

\newcommand{\Lime}{{\color{lime!60!black} lime\,}}

\newcommand{\rDim}{\textup{{\color{red!60!black}red}-imension }}
\newcommand{\lDim}{\textup{{\color{lime!60!black}limed}-imension }}

\begin{document}

\title{Linear Algebra In A Jiffy}
\markright{Linear Algebra In A Jiffy}
\author{Leo Livshits}  
\maketitle

\begin{center}
 {\emph {Dedicated to Heydar Radjavi on the occasion of his $89$-th birthday.}}
\end{center}

\begin{abstract}
Every $n$-tuple in $\F^{n}$ has a first non-zero entry and a last non-zero entry. What do the positions of such entries in the elements of a subspace \Fam{W} of $\F^{n}$ reveal about \Fam{W}? It turns out, a great deal!
This insight offers a beeline to the fundamental results on general bases and dimension (including dimensions of complements, rank of a transpose, Rank-Nullity Theorem, Full Rank Factorization, and existence/uniqueness of RREF), bypassing some of the traditional stumbling blocks and time sinks. \end{abstract}

\section{Introduction}
This article is about doing old things in a new way. It is aimed at mathematicians who teach linear algebra, and it is laid out in a fairly formal way, for efficiency sake.\footnote{It is a useful feature of the arguments that they are short and are easily converted into a concrete and visual form that students find pleasantly comprehensible and convincing, despite the absence of full generality and rigorous formalism. This is how the material is presented in our forthcoming textbook ``Not Your Grandpa's First Introduction To Linear Algebra.''} The reader who decides to embark on this trip is advised to exercise a touch of amnesia; otherwise the tenets of the standard approach will not let go of the mind.

Modern developments suggest that initial focus in introductory linear algebra should be placed on $\R^{n}$ and $\C^{n}$, with concepts of orthogonality and least squares playing an earlier and more central role, and Singular Value Decomposition done justice within the first course. To be able to implement this without skipping arguments that justify claims, or replacing proper development with commandments, something has to go overboard. 

Our approach establishes basic fundamental results quickly, bypassing traditional stumbling blocks and time sinks. 
The gained efficiency is due to an early introduction of a combinatorial concept of right/left dimension of a subspace, requiring no prior exposure to linear independence, linear systems or row-reduction.

In a course, the material that should precede our point of departure, would include an introduction to the operations on $\R^{n}$ (or $\C^{n}$, or $\F^{n}$), the concept of a subspace, and the fact that spans and ortho-complements generate subspaces. No  more than that.


 The article is organized as follows. In section 2 we introduce right/left-standard indices. This is a new concept that serves as a trailhead for the shortcut. Section 3 deals with the right/left-basic elements of subspaces and canonical coordinate systems. Section 4 treats dimensionality.  In section 5 we connect the canonical coordinate systems of a subspace with those of its ortho-complement.

The material in sections 6 and 7 assumes some rudimentary familiarity with matrix algebra, including transposition. Again, only the basics are required. In these sections we apply the new methods to establish the equality of row and column ranks, Rank-Nullity theorem, full rank factorization, and existence/uniqueness of the RREF, all as immediate by-products. All of this happens quickly and at a low cost.

Section 8 describes a recursive form of Gauss-Jordan elimination reminiscent of Gram-Schmidt process, that can be used to construct the RREF of a matrix.

Section 9 of the article deals with a combinatorial question that is not usually a part of the course. If one were to mark all of the positions where the elements of a subspace have their first non-zero entry and their last non-zero entry, what configurations can be thus produced?
%

Throughout this article, $\F$ represents an arbitrary field. We interpret the vectors of $\F^{n}$ as lists, and as such these can be presented vertically, horizontally, or even diagonally, as long as the ordering is clear. The entries of a vector in $\F^{n}$ appear in its $n$ {\bf \textit{positions}}, indexed $1$ through $n$. The $i$-th entry of $X\in\F^{n}$ is the element of $\F$ that appears in the $i$-th position of $X$, and it is denoted by ${\BM{X(i)}}$. The \textit{\textbf{$\ensuremath{\BM{\it{k}}}$-th standard basis vector}} $\ensuremath{\BM{\it{E_{k}}}}$ in $\F^{n}$ is the list that has $1$ as its $k$-th entry and zeros in all other positions. 

The last non-zero entry of a vector in $\F^{n}$ is said to be the \textit{\textbf{terminating entry}} of that vector. Obviously the zero vector has none such. If the terminating entry of $X\in\F^{n}$ occurs in the $j$-th position, we say that $X$ \textit{\textbf{terminates}} in the $j$-th position. Similarly the first non-zero entry of an element of $\F^{n}$ is its \textit{\textbf{originating entry}}, and $X$ \textit{\textbf{originates}} in the $j$-th position if its originating entry occurs in that position.

%

\section{Right/Left-Standard Indices}\mbox{}
Let \Fam{W} be a subspace of $\F^{n}$. If \Fam{W} contains an element that terminates in the $j$-th position, we say that $j$ is a \textit{{right-standard index}} of \Fam{W}, and the $j$-th position is a \textit{{right-standard position}} for \Fam{W}. For example, if \Fam{W} contains an element with no zero entries, then the $n$-th position is a right-standard position for \Fam{W}. 
Since subspaces are closed under scaling, a position is right-standard for $\Fam{W}$ if and only if $\Fam{W}$ contains an element which terminates with a $1$ in that position. 

To avoid tongue-twisting, we will trade in the term ``right-standard'' for ``{\textbf{\color{red!60!black} red}}''. Why {\color{red!60!black} red}? Well, both terms start with an ${\bf r}$, ``{\color{red!60!black} red}'' is a short word, and \Red entries are bound to leave \Red residue, marking the \Red positions.

Similarly, if \Fam{W} contains an element that originates in the $j$-th position, then $j$ is a \textit{{left-standard index}} of \Fam{W}, and the $j$-th position is a \textit{{left-standard position}} for \Fam{W}. We write ``{\color{lime!50!black} lime}'' for ``left-standard''.

Since results involving {\color{lime!50!black} lime} indices, entries and positions are mostly analogous to those  involving the \Red versions thereof, we will just focus on the latter, unless we specify otherwise.

The number of {\color{red!60!black} red} positions for \Fam{W} is said to be the  {\bf \textit{{\color{red!60!black}red}-imension}} of \Fam{W}. Obviously the \rDim of the trivial subspace $\{O_{n}\}$ is zero. We choose this (temporary) strange terminology to highlight the fact that we are not cheating by somehow sneaking in the usual notion of dimension before the underlying concepts and results are established.
 

Since any non-zero element of \Fam{W} terminates in a {\color{red!60!black} red} position, 
\MD{\small the only element of \Fam{W} that has zeros in all {\color{red!60!black} red} positions of \Fam{W} is the zero vector.\hfill \textup{(*)}} 
As \Fam{W} is closed under subtraction,
it follows that \MD{\small two elements of \Fam{W} are equal exactly when they coincide in all {\color{red!60!black} red} positions of \Fam{W}.\ {\footnotesize\textup{(**)}}}

In particular, if $i$ is {\emph {the smallest}} {\color{red!60!black} red} index of \Fam{W}, then there is exactly one element in $\Fam{W}$ which terminates with a $1$ in the $i$-th position.

The existence of the \Red indices for a non-trivial subspace of $\F^{n}$ is undeniable, {\emph {no matter how the subspace is defined}}. Finding what these indices are in a particular situation is a different story, as is often the case in linear algebra. This does not diminish the value of the approach, which is designed to beeline through the introductory fundamentals of linear algebra in an honest way.

\section{Red-Basic Elements}\mbox{}
\begin{Thm}[{\color{red!60!black} Red}-Basic Elements]\label{T3}
When \Fam{W} is a non-trivial subspace of $\F^{n}$, for each {\color{red!60!black} red} index $i$ of \Fam{W} there exists a {{unique}} element $W_{i}\in\Fam{W}$ that terminates with a $1$ in the $i$-th position, and has $0$'s in all preceding {\color{red!60!black} red} positions of \Fam{W}.
\end{Thm}

\begin{proof}
Let us start with the uniqueness, as this is an easier task. If $X$ and $Y$ were distinct elements of \Fam{W} with the described property for a {\color{red!60!black} red} index $i$, their difference would be a non-zero element of \Fam{W} that has zeros in all {\color{red!60!black} red} positions of \Fam{W}, which is not possible, by (**).

On to the existence. As we have already observed, the claim holds true for the smallest \Red index $i_{_{1}}$ of \Fam{W}. We proceed inductively. 
Suppose that $W_{_{i_{1}}}, W_{_{i_{2}}}, \ldots , W_{_{i_{k}}}$ have been constructed for the smallest $k$ \Red indices $i_{_{1}}, i_{_{2}},  \ldots, i_{_{k}}$ of \Fam{W}, and $i_{_{k+1}}$ is the next \Red index of \Fam{W}. 

\Fam{W} contains an element $X$ which terminates with a $1$ in the $i_{_{k+1}}$\!-st position. It follows that   
\begin{equation*}
X-X\left(i_{1}\right)W_{_{i_{1}}}-X\left(i_{2}\right)W_{_{i_{2}}}-\cdots-X\left(i_{k}\right)W_{_{i_{k}}}
\end{equation*}
is an element of \Fam{W} that has zero entries in positions indexed by  $i_{_{1}}, i_{_{2}},  \ldots, i_{_{k}}$, and terminates with a $1$ in the $i_{_{k+1}}$\!-st position. This is the $W_{{i_{k+1}}}$ we seek. 
\end{proof}

For each \Red index $i$ of \Fam{W}, we will refer to the element $\BM{W_{i}}$ described in theorem \ref{T3} as {\textit{\textbf  {the {\color{red!60!black} red}-basic element}}}  of \Fam{W} corresponding to $i$.  
Restricted to the {\color{red!60!black} red} positions of \Fam{W} the {\color{red!60!black} red}-basic elements of \Fam{W} look like the standard basis $m$-tuples $E_{1}, E_{2},  \ldots, E_{m}$, where $m$ is the \rDim of \Fam{W}.

\newpage

\begin{Thm}[{\color{red!60!black} Red} Entries Are Free]\label{T1;31-1-24}

Suppose that a subspace \Fam{W} of $\F^{n}$ has  \rDim $k$, and $a_{1}, a_{2},  \ldots, a_{k}\in\F$. Then there exists {exactly one} element $X\in\Fam{W}$ with $a_{1}, a_{2},  \ldots, a_{k}$ as its entries in the {\color{red!60!black} red} positions of \Fam{W}  (in that order). 
\end{Thm}

\begin{proof}
We have already noted in (**) that there cannot be two distinct elements of \Fam{W} with the described property. All we need to show is that there is always at least one such.

This is not difficult, since on the {\color{red!60!black} red} positions of \Fam{W} the {\color{red!60!black} red}-basic elements of \Fam{W} look like the standard basis vectors $E_{1}, E_{2},  \ldots, E_{k}$ of $\F^{k}$.

If $i_{_{1}}< i_{_{2}}<  \cdots< i_{_{k}}$ are the \Red indices of \Fam{W}, then the element 
\begin{equation*}
X= a_{1}W_{_{i_{1}}}+a_{2}W_{_{i_{2}}}+\cdots+a_{k}W_{_{i_{k}}}\in\Fam{W}
\end{equation*}
satisfies the required condition of having $a_{1}, a_{2},  \ldots, a_{k}$ as its entries in the {\color{red!60!black} red} positions of \Fam{W}.
\end{proof}

\begin{Dfn}[Coordinate Systems]\label{d1;24-2-24}
We say that a list $X_{1}, X_{2},  \ldots, X_{k}$ in a subspace \Fam{W} of $\F^{n}$ generates a {\bf \textit{coordinate system}} for \Fam{W} if every element of \Fam{W} can be expressed uniquely as a linear combination of $X_{1}, X_{2},  \ldots, X_{k}$.
\end{Dfn}

\begin{Cor}[{\color{red!60!black} Red}-basic elements generate a coordinate system]\label{c3}
If \Fam{W} is a non-trivial subspace of $\F^{n}$ then the {\color{red!60!black} red}-basic elements of \Fam{W} generate a coordinate system of \Fam{W}.
\end{Cor}

\begin{proof}
We use the fact that on the {\color{red!60!black} red} positions of \Fam{W} the {\color{red!60!black} red}-basic elements look like the standard basis vectors. 

If $i_{_{1}}< i_{_{2}},<  \cdots< i_{_{k}}$ are the \Red indices of \Fam{W}, and $X\in\Fam{W}$ then  
\begin{equation*}
X\left(i_{1}\right)W_{_{i_{1}}}+X\left(i_{2}\right)W_{_{i_{2}}}+\cdots+X\left(i_{k}\right)W_{_{i_{k}}}
\end{equation*}
is the unique element of \Fam{W} with $X\left(i_{1}\right), X\left(i_{2}\right), \ldots, X\left(i_{k}\right)$ as its entries in the \Red positions of \Fam{W}. Of course $X$ itself has these same entries in those same positions. By (**),
 \begin{equation*}
X\left(i_{1}\right)W_{_{i_{1}}}+X\left(i_{2}\right)W_{_{i_{2}}}+\cdots+X\left(i_{k}\right)W_{_{i_{k}}}
\end{equation*}
is the one and only linear combination of the \Red\!-basic elements of \Fam{W} that equals $X$.
\end{proof}

\MD{\small To express an element $X$ in \Fam{W} as a linear combination of the {\color{red!60!black} red}-basic elements of \Fam{W}, one uses the entries of $X$ in the {\color{red!60!black} red}  positions of \Fam{W}  as coefficients.}

In particular, if   $W_{i_{1}}, W_{i_{2}},  \ldots, W_{i_{k}}$ are the {\Red\!-basic} elements of \Fam{W}, testing whether a given $n$-tuple $X$ is in \Fam{W} is as easy as pie: 
simply test the equality
\begin{equation}\label{e1}
X=X\left(i_{1}\right)W_{i_{1}}+\cdots+X\left(i_{k}\right)W_{i_{k}}.
\end{equation}


\section{Red Indices Give Dimension}\label{sec4;19-2-24}\mbox{}
The following result is standard and has an easy enough standard proof, which we omit for the sake of brevity.

\begin{Thm}[Coordinate system criterion]\label{t1;18-2-24}
A list $X_{1}, X_{2},  \ldots, X_{k}$ in a subspace \Fam{W} of $\F^{n}$ generates a coordinate system for \Fam{W} (see  \ref{d1;24-2-24}) exactly when it spans \Fam{W}, does not start with a zero vector, and has no elements that are linear combinations of the preceding elements. 
\end{Thm}

In view of theorem \ref{t1;18-2-24}, we can use the traditional term ``{\bf \textit{basis}}'' and  our term ``coordinate system'' interchangeably. The {\bf  \textit{{\color{red!60!black} red} basis}} of \Fam{W} is the coordinate system formed by its \Red\!-basic elements listed in the order of the corresponding \Red indices (see corollary \ref{c3}).

\begin{Thm}[Monotonicity of \redim]\label{C5}

If $\Fam{W} \preceq \Fam{V}\preceq \F^{n}$, then \ $\rDim\!\!\left(\Fam{W}\right)\leq \rDim\!\!\left(\Fam{V}\right),$
and the equality of {\color{red!60!black} red}-imensions holds only when $\Fam{W}= \Fam{V}$.
\end{Thm}

\begin{proof}
The {\color{red!60!black} red} indices of a subspace are determined by the terminal positions of its elements. Through enlarging the subspace one can only add to the collection of such indices.

If no new {\color{red!60!black} red} indices have been added in the process, then the {\color{red!60!black} red}-basic elements of \Fam{W} must be exactly the {\color{red!60!black} red}-basic elements of \Fam{V}, which means that \Fam{W} and \Fam{V} have the same {\color{red!60!black} red} basis and are therefore equal.
\end{proof}

\begin{Thm}[Step-up bound for \redim]\label{T2}
For any $X_{1}, X_{2},  \ldots\!, X_{k}, Y\in\F^{n}$,
\begin{multline*}
\rDim\!\!\Big(\Span{X_{1}, X_{2},  \ldots, X_{k}, Y}\!\!\Big)\\ \leq \rDim\!\!\Big(\Span{X_{1}, X_{2},  \ldots, X_{k}}\!\!\Big)\ +1.
\end{multline*}
\end{Thm}

\begin{proof}
Let us denote the span of $X_{1}, X_{2},  \ldots, X_{k}$ by \Fam{W}, and the span of \\ $X_{1}, X_{2},  \ldots, X_{k}, Y$ by \Fam{V}. Elements of \Fam{V} have the form $Z+\alpha Y$, where $Z\in\Fam{W}$, $\alpha\in \F$.

We need to demonstrate that we could not have created two new {\color{red!60!black} red} indices for the span, by attaching $Y$ to the original list $X_{1}, X_{2},  \ldots, X_{k}$. Had we done so, there would be elements $Z_{1}+\alpha Y$ and $Z_{2}+\beta Y$ that terminate in distinct non-{\color{red!60!black} red} positions for \Fam{W} (say, $i$-th and $j$-th positions, with $i<j$). 

Neither $\alpha$ nor $\beta$ can be zero. Consequently $\beta Z_{1}+\alpha\beta Y$ and $\alpha Z_{2}+\alpha\beta Y$ still terminate in the $i$-th and $j$-th positions respectively. Therefore their difference, which is $\alpha Z_{2}-\beta Z_{1}$, terminates in the $j$-th position (since $i<j$). Yet this difference is an element of \Fam{W}, and hence cannot terminate in a non-{\color{red!60!black} red} position of \Fam{W}. This contradiction completes the proof.
\end{proof}

\begin{Thm}[\Redim is the dimension]\label{T1;29-2-24}
Every basis of \Fam{W} has exactly the same number of elements as the {\color{red!60!black} red} basis of \Fam{W}.\footnote{It is worth reiterating that naturally we have made neither use of nor assumption about the equality of cardinalities of various bases of a given subspace. The equality is an easy byproduct of our tack.}\BSK
\end{Thm}

\begin{proof}
By theorem \ref{t1;18-2-24}, a list $X_{1}, X_{2},  \ldots, X_{m}$ in a non-trivial subspace $\Fam{W}$ of $\F^{n}$ is a basis of \Fam{W} exactly when 
\begin{multline*}
\{O\}\subsetneq \Span{X_{1}}
\subsetneq\Span{X_{1}, X_{2}}\subsetneq\cdots\subsetneq\Span{X_{1}, X_{2},  \ldots, X_{m}}=\Fam{W}.
\end{multline*}
By theorem \ref{C5}, this, apart from the requirement that $X_{i}\in\Fam{W}$, can be restated in terms of \rDim:
\begin{multline*}
1= \rDim\!\!\Big(\Span{X_{1}}\!\!\Big)<\rDim\!\!\Big(\Span{X_{1}, X_{2}}\!\!\Big)<\cdots\\ <\rDim\!\!\Big(\Span{X_{1}, X_{2},  \ldots, X_{m}}\!\!\Big)=\rDim\!\!\left(\Fam{W}\right).
\end{multline*}
In such a case each step-up in \rDim must be exactly $1$, by theorem \ref{T2}, and the proof is complete.
\end{proof}

\begin{Cor}
Every subspace \Fam{W} of $\F^{n}$ has as many \Lime indices and \Red ones.\footnote{Of course the reader will have noticed that the linear isomorphism on $\F^{n}$ that reverses the order of the entries in each $n$-tuple, maps a subspace \Fam{W} to a subspace $\Fam{W}^{\REV}$ in such a way that the {\color{red!60!black} red} indices of \Fam{W} become the {\color{lime!50!black} lime} indices of $\Fam{W}^{\REV}$. This demonstrates that  \rDim and \lDim of a subspace are equal, but obviously by  ``borrowing from the future''. 

A compelling reason for introducing {\color{red!60!black} red}/{\color{lime!50!black} lime} indices is to provide an efficient alternate path for the first linear algebra course. The {\color{red!60!black} red}/{\color{lime!50!black} lime} indices are designed to be encountered well before linear functions make their official appearance, or the traditional concept of dimension is introduced. 
}\BSK
\end{Cor}



Another easy observation is the following:
\begin{Obs}\LL{Obs1;28-5-25}
If a list that spans a subspace \Fam{W} of $\F^{n}$ looks like a \Red basis then it is the \Red basis of \Fam{W}.

What we mean by ``looks like a \Red basis'' is this:
\begin{enumerate}
\item Every vector on the list originates with a $1$.
\I Every vector on the list has a zero entry in every position where another vector on the list has its originating $1$.
\end{enumerate}
\end{Obs}

%
%
%

\section{Red/Lime Bases And $\bullet$-Complements}

We will write $X\bullet Y$ to indicate the standard symmetric bilinear form, commonly known as the dot product, for elements of $\F^{n}$. For each collection \Fam{S} in $\F^{n}$, we define its \textit{\textbf{$\bullet$-complement}} to be
\begin{equation*}
\Fam{S}^{^{\bullet}}=\Set{X\in\F^{n}}{\forall\, Y\in\Fam{S}: X\bullet Y=0}.
\end{equation*}
When $\F=\R$, $\bullet$-complement captures the geometry of perpendicularity. \footnote{For the standard inner product on $\C^{n}$ see remark \ref{R1}.}

\begin{Thm}[Duality of {\color{red!60!black} red}/{\color{lime!60!black} lime} indices of $\bullet$-complements]\label{T1;18-1-24}

For any subspace $\Fam{W}$ of $\F^{n}$, the {\color{lime!50!black} lime} indices of $\Fam{W}^{^{\bullet}}$ are exactly the non-{\color{red!60!black} red} indices of \Fam{W}.
\end{Thm}

\begin{proof}
If $i$ is a {\color{red!60!black} red} index of \Fam{W}, it cannot be a {\color{lime!50!black} lime} index of $\Fam{W}^{^{\bullet}}$, since otherwise \Fam{W} would contain an element that terminates with a $1$ in the $i$-th position, while $\Fam{W}^{^{\bullet}}$ would contain an element that originates with a $1$ in the $i$-th position. The dot product of these two elements would be $1$, contrary the definition of $\Fam{W}^{^{\bullet}}$.
%
%
%
%
%
%

If $i_{{o}}$ is a non-{\color{red!60!black} red} index of \Fam{W}, then there are two scenarios to consider. In the first scenario there are no {\color{red!60!black} red} indices of \Fam{W} that are greater than $i_{{o}}$. In that case every element of  \Fam{W} terminates before the $i_{{o}}$-th position. Therefore the standard basis vector $E_{i_{{o}}}$ is in $\Fam{W}^{^{\bullet}}$, which shows that $i_{{o}}$ is a {\color{lime!50!black} lime} index of $\Fam{W}^{^{\bullet}}$.

In the remaining scenario $i_{o}$ is followed by some {\color{red!60!black} red} indices, say $i_{_{1}}, i_{_{2}},  \ldots, i_{_{k}}$. Let us write $a_{_{j}}$ for the $i_{o}$-th entry 
 of $W_{i_{j}}$. 

We shall construct an element $Z$ of $\Fam{W}^{^{\bullet}}$ which originates in the $i_{o}$-th position, so that $i_{o}$ is a {\color{lime!50!black} lime} index of $\Fam{W}^{^{\bullet}}$.\footnote{This proof is followed by an illustration of the construction.} The design of $Z$ is such that for each {\color{red!60!black} red}-basic element $W_{i}$ of \Fam{W}, there are at most two positions where both $W_{i}$ and $Z$ have non-zero entries.

We let
\begin{equation*}
Z(m)=\left\{\begin{array}{cl}
1& \text{if } m=i_{o}\\
-a_{_{j}}& \text{if } m=i_{j}; 1\leq j\leq k\\
0 & \text{otherwise}
\end{array}\right..
\end{equation*}
It is easy to see that $Z\bullet W=0$ for every {\color{red!60!black} red}-basic element $W_{i}$ of \Fam{W}, and since these form a basis of \Fam{W}, $Z\in \Fam{W}^{^{\bullet}}$.
\end{proof}

Here is an illustration of the construction in the proof of theorem \ref{T1;18-1-24}, when $i_{o}=2$ and $i_{_{1}}, i_{_{2}},  \ldots, i_{_{k}}$ are $5,8,9,13,\ldots,p$:\Vm{-0.5}
\begin{equation*}\label{}
\resizebox{\textwidth}{!}{$\begin{array}{rcccccccccccccccccc}
W_{5}&=&(\star,& {\BM{a_{5}}},& \star,& \star,& {\color{red!60!black} 1},& 0,& 0,&{\color{red!60!black} 0},& {\color{red!60!black} 0},& 0,& 0,& 0,& {\color{red!60!black} 0},& \ldots,& {\color{red!60!black} 0},& \ldots, & 0)\\
W_{8}&=&(\star,& {\BM{a_{8}}},& \star,& \star,& {\color{red!60!black} 0},& \star,& \star,& {\color{red!60!black} 1},& {\color{red!60!black} 0},& 0,& 0,& 0,& {\color{red!60!black} 0},& \ldots,& {\color{red!60!black} 0},& \ldots &0)\\
W_{9}&=&(\star,& \BM{a_{9}},& \star,& \star,& {\color{red!60!black} 0},& \star,& \star,& {\color{red!60!black} 0},& {\color{red!60!black} 1},& 0,& 0,& 0,& {\color{red!60!black} 0},& \ldots,& {\color{red!60!black} 0},& \ldots & 0)\\ : \\
W_{p}&=&(\star,& \BM{a_{p}},& \star,& \star,& {\color{red!60!black} 0},& \star,& \star,& {\color{red!60!black} 0},& {\color{red!60!black} 0},& \star,& \star,& \star,& {\color{red!60!black} 0},& \ldots,& {\color{red!60!black} 1},& \ldots & 0)
\\ \\  
Z&=&(0,& {\BM{{\color{lime!60!black} 1}}},& 0,& 0,& {\BM{-a_{5}}},& 0,& 0,& {\BM{-a_{8}}},& {\BM{-a_{9}}},& 0,& 0,& 0,& {\BM{-a_{13}}},& \ldots,& {\BM{-a_{p}}},& \ldots & 0)\\
&&&{\color{lime!60!black} \uparrow} &&&{\color{red!60!black} \uparrow}&&&{\color{red!60!black} \uparrow}&{\color{red!60!black} \uparrow}&&&&{\color{red!60!black} \uparrow}&&{\color{red!60!black} \uparrow}&&
\end{array}$}
\end{equation*}

\begin{Obs}\label{r1;19-1-24}
It is worth noting that the element $Z$ constructed in the proof of theorem \ref{T1;18-1-24} is exactly the {\color{lime!50!black} lime}-basic element of $\Fam{W}^{^{\bullet}}$ corresponding to the {\color{lime!50!black} lime} index $i_{o}$ of $\Fam{W}^{^{\bullet}}$. This shows that\MD{\small the {\color{lime!50!black} lime} basis of $\Fam{W}^{^{\bullet}}$ can be simply read off from the {\color{red!60!black} red} basis of \Fam{W}, and vice versa. }
\end{Obs}

\begin{Cor}[Dimensions of $\bullet$-complements]\label{C1;16-2-24}
For any subspace \Fam{W} of $\F^{n}$, $\Dim{\Fam{W}}+\Dim{\Fam{W}^{^{\bullet}}}=n$, and consequently\\ $\left(\Fam{W}^{^{\bullet}}\right)^{^{\!\!\!\bullet}}=\Fam{W}.$
\end{Cor}

\begin{proof}
We have verified that the {\color{lime!50!black} lime} indices of $\Fam{W}^{^{\bullet}}$ are exactly the non-{\color{red!60!black} red} indices of \Fam{W}. Of course this shows that $\Dim{\Fam{W}^{^{\bullet}}}=n-\Dim{\Fam{W}}.$

To establish the remaining claim, note that $\Fam{W}$ is naturally included in $\left(\Fam{W}^{^{\bullet}}\right)^{^{\!\!\!\bullet}}$, and the two subspaces have the same dimensions.
\end{proof}

\begin{Obs}\label{R1}
An analogous argument (cf. proof of theorem \ref{T1;18-1-24}) demonstrates similar results in the case $\F=\C$, when $\bullet$ is replaced by the standard inner product. This can be viewed as a particular case of a more general claim. If $\varphi:\F\lra\F$ is an involutive isomorphism, and we define $\Fam{P}:\F^{n}\times\F^{n}\lra\F$ by
\begin{equation*}
\Fam{P}(X,Y)=X\bullet \varphi^{(n)}(Y),
\end{equation*}
where
\begin{equation*}
\varphi^{(n)}{\scriptstyle{
\left(\begin{array}{c}
\scriptstyle{y_{1}}\\
\scriptstyle{y_{2}}\\
: \\
\scriptstyle{y_{n}}
\end{array}\right)}}
={\scriptstyle{
\left(\begin{array}{c}
\scriptstyle{\varphi}\left(\scriptstyle{y_{1}}\right)\\
\scriptstyle{\varphi}\left(\scriptstyle{y_{2}}\right)\\
: \\
\scriptstyle{\varphi}\left(\scriptstyle{y_{n}}\right)
\end{array}\right)}},
\end{equation*}
then analogous conclusions can be drawn regarding $\Fam{W}^{^{\Fam{P}}}$ and $\left(\Fam{W}^{^{\Fam{P}}}\right)^{^{\!\!\!\Fam{P}}}$.
\end{Obs}

\section{Red/Lime Bases And Ranges/Nullspaces}\mbox{}
Let us begin with a brief review of some concepts and terminology pertaining to matrices, which would have to be introduced prior to the material we are about to discuss. 

When one applies an $n\times m$ matrix  $\BMR{A}$ to an $m$-tuple $X$, the resulting $n$-tuple $\BMR{A}(X)$ can be described in terms of the rows $R_{1}, R_{2},  \ldots, R_{n}$ of $\BMR{A}$, as well as in terms of the columns $C_{1}, C_{2},  \ldots, C_{m}$ of $\BMR{A}$:
\begin{enumerate}
\item the $i$-th entry of $\BMR{A}(X)$ is the dot products $R_{i}\bullet X$;\SSK
\item $\BMR{A}(X)=X(1)C_{1}+X(2)C_{2}+\cdots+X(m)C_{m}$.
\end{enumerate}
For convenience let us refer to these as a {\bf row-centric} and a {\bf column-centric} descriptions of $\BMR{A}(X)$, respectively.

The column-centric perspective shows that, interpreted as a linear function acting from $\F^{m}$ to $\F^{n}$, $\BMR{A}$ is invertible (i.e. bijective) if and only if every element of $\F^{n}$ can be expressed uniquely as a linear combination of the columns of $\BMR{A}$; in other words, if and only if the columns of $\BMR{A}$ form basis of $\F^{n}$ (see definition \ref{d1;24-2-24} and theorem \ref{t1;18-2-24}). As we have already established in theorem \ref{T1;29-2-24}, every basis of $\F^{n}$ has $n$ elements. Hence invertible matrices are square.

The row-centric perspective shows that $\BMR{A}(X)$ is a zero vector exactly when $X$ is in the $\bullet$-complement  of rows of $\BMR{A}$ within $\F^{m}$; in other words, when $X$ is in the $\bullet$-complement of the span of these rows, known as the {\bf row space} of $\BMR{A}$. This states that \MD{\small $\text{Nullspace}{\left(\BMR{A}\right)}$ and $\text{Row Space}(\BMR{A})$ are $\bullet$-complements within $\F^{m}$.}

While we are on the topic, let us remind the reader that the span of the columns of a matrix $\BMR{A}$ is said to be its {\bf column space}, and is exactly the range of $\BMR{A}$, when $\BMR{A}$ is interpreted as a linear function from $\F^{m}$ to $\F^{n}$. {\bf Column rank} of $\BMR{A}$ is defined to be the dimension of the column space of $\BMR{A}$, and {\bf row rank} is the dimension of the row space. Since columns of $\BMR{A}^{T}$ (the transpose of $\BMR{A}$) are rows of $\BMR{A}$, column space of $\BMR{A}^{T}$ equals row space of $\BMR{A}$, and row rank of $\BMR{A}$ equals column rank of $\BMR{A}^{T}$.

By corollary \ref{C1;16-2-24},\Vm{-1}
\begin{multline}\label{e1;25-2-24}
m=\Nlty{\BMR{A}}+\RRank{\BMR{A}}=\Nlty{\BMR{A}}+\CRank{\BMR{A}^{T}}.
\end{multline}

Using the column-centric perspective and theorem \ref{T3}, it is not hard to see that each of the following claims regarding an index $i$ and a matrix $\BMR{A}$ with columns $C_{_{1}}, C_{_{2}},  \ldots, C_{_{m}}$ is equivalent to the rest.
\begin{enumerate}
\I $C_{i}$ is a linear combination of the preceding columns of $\BMR{A}$.
\I There are scalars $a_{1}, a_{2},  \ldots, a_{i-1}\in\F$ such that
\begin{equation*}
O=a_{1}C_{1}+a_{2}C_{2}+\cdots+a_{i-1}C_{i-1}+1C_{i}
\end{equation*}
\I There is a vector $X$ in the nullspace of $\BMR{A}$ that terminates with a $1$ in the $i$-th position.
\item $i$ is a \Red index for the nullspace of $\BMR{A}$.
\I There is a vector $W_{i}$ in the nullspace of $\BMR{A}$ that terminates with a $1$ in the $i$-th position, and has zeros in all of the preceding \Red positions.
\I $C_{i}$ is a linear combination of the preceding columns of $\BMR{A}$ that correspond to non-\Red indices of the nullspace of $\BMR{A}$.
\end{enumerate}
Since $\text{Nullspace}{\left(\BMR{A}\right)}$ and $\text{Row Space}(\BMR{A})$ are $\bullet$-complements, the \Red indices of the nullspace of $\BMR{A}$ are exactly the non-\Lime indices of the row space of $\BMR{A}$ (by theorem \ref{T1;18-1-24}). It follows, by theorem \ref{t1;18-2-24}, that
 \MD{\small if $i_{1}< i_{2}<  \cdots< i_{k}$ are the {\color{lime!50!black} lime} (resp. \Red\!) indices of the row space of $\BMR{A}$,\hfill\mbox{}\\ then the columns $C_{i_{1}}, C_{i_{2}}, \ldots, C_{i_{k}}$ form a basis of the column space of $\BMR{A}$. \footnote{Of course the roles of rows and columns can be reversed here  through transposition.}\hfill 
\textup{($\dagger$)}\BSK}
Therefore
 \begin{equation*}
\RRank{\BMR{A}}=\CRank{\BMR{A}},
\end{equation*}
and we just use $\Rank{\BMR{A}}$ for both.
In particular,
\begin{equation}\label{e3;19-1-24}
\Rank{\BMR{A}^{T}}=\Rank{\BMR{A}}.
\end{equation} 
By \eqref{e1;25-2-24}, we have arrived that
\begin{equation*}
\Rank{\BMR{A}}+\Nlty{\BMR{A}}=m  \ \ \ \ \text{(Rank-Nullity Formula).}
\end{equation*}

\section{Red/Lime Bases and RREF/RCEF}

Next we discuss the connection between RREF/RCEF and {\color{lime!50!black} lime} bases, touching base with Grinberg's ``gauche perspective'' \cite{Grinberg}. Existence and uniqueness of $\text{RREF}(\BMR{A})$, for any matrix $\BMR{A}$, is an immediate by-product of our approach.

Our path to RREF/RCEF passes through the Full Rank factorization, which is fundamental in its own right. If the reader has not encountered it previously, we recommend treating it as a generalization of the fact that matrices of rank $1$ are exactly those of the form $[X]_{_{n\times 1}}[Y]^{^{T}}_{_{m\times 1}}$.

\begin{Thm}[Full Rank factorization using rows of $\BMR{A}$]\label{T2;18-1-24}
Suppose that a matrix $\BMR{A}\in\Mn{n\times m}$ has rank $r>0$, and $W_{i_{1}}, W_{i_{2}}, \ldots , W_{i_{r}}$ is the \Lime basis of the columns space of $\BMR{A}$. 

Let $\BMR{B}$  be the matrix whose columns are $W_{i_{1}}, W_{i_{2}}, \ldots , W_{i_{r}}$, and let $\BMR{G}$ be the matrix whose rows are the rows $R_{i_{1}}, R_{i_{2}}, \ldots , R_{i_{r}}$ of $\BMR{A}$. \footnote{By ($\dagger$) these rows form a basis of the row space of $\BMR{A}$.}
Then 
\begin{equation*}
\BMR{A}=\BMR{BG} \ \ \AND\ \ \Rank{\BMR{B}}=r=\Rank{\BMR{G}}.
\end{equation*}
\smallskip
\end{Thm}

\begin{proof}
The rank condition is immediate for $\BMR{B}$, and for $\BMR{G}$ it follows from ($\dagger$).

The fact that $\BMR{A}=\BMR{BG}$ is best seen one column at a time, via the column-centric perspective. Each column $C_{j}$ of $\BMR{A}$ is a linear combination of $W_{i_{1}}, W_{i_{2}}, \ldots , W_{i_{r}}$, i.e. of the columns of $\BMR{B}$, with coefficients being the entries of $C_{j}$ corresponding to the indices $i_{1}< i_{2}<  \cdots< i_{r}$\ \ (see \eqref{e1}\!\! ). These entries in the columns $C_{1}, C_{2},  \ldots, C_{m}$ form rows of $\BMR{A}$ that correspond to the indices $i_{1}< i_{2}<  \cdots< i_{r}$; i.e. the rows of $\BMR{G}$.
\end{proof}

A matrix whose rows are the {\color{lime!50!black} lime} basis of its row space, in the increasing order of indices, followed perhaps by some zero rows at the bottom, is said to be in {\bf {Reduced Row Echelon Form}} (RREF). 

If one switches ``rows'' to ``columns'' (and ``bottom'' to ``right'') in the paragraph above, one arrives at {\bf {Reduced Column Echelon Form}} (RCEF).

It is not hard to see that the transpose of an RREF matrix is an RCEF matrix, and conversely.

The use {\color{lime!50!black} lime} bases renders the question of uniqueness of RREF/RCEF trivial: 
for a given $m\times n$ matrix $\BMR{A}$, we write $\textbf{RREF}(\BMR{A})$ for the unique  $m\times n$ RREF matrix with the same row space as $\BMR{A}$. 

Since row space of a matrix is the $\bullet$-complement of its nullspace,  $\RREF{\BMR{A}}$ is the unique  $m\times n$ RREF matrix with the same nullspace as $\BMR{A}$.\footnote{In particular, any RREF matrix is completely determined by its nullspace.} By remark \ref{r1;19-1-24}, one can read off the {\color{red!60!black} red} basis of the nullspace of an RREF matrix from the matrix itself with minimal effort.

Similarly, we write $\textbf{RCEF}(\BMR{A})$ for the unique  $m\times n$ RCEF matrix with the same column space as $\BMR{A}$. Since column space of \BMR{A} equals row space of  $\BMR{A}^{T}$, the two spaces share the \Lime basis, and it follows that
\begin{equation}\label{e2;19-1-24}
\Big(\RCEF{\BMR{A}}\!\!\Big)^{\!T}=\RREF{\BMR{A}^{T}}.
\end{equation}

\begin{Cor}[RCEF factorization]\label{c1;19-1-24}
Suppose that $\BMR{A}$ is an $n\times m$ non-zero matrix, and $i_{1}< i_{2}<  \cdots< i_{r}$ are the \Lime indices of the column space of $\BMR{A}$. Then 
\begin{equation}\label{e1;5-3-24}
\BMR{A}=\RCEF{\BMR{A}}\circ \BMR{S},
\end{equation}
where $\BMR{S}$ is the matrix whose rows are the rows $R_{i_{1}}, R_{i_{2}}, \ldots , R_{i_{r}}$ of $\BMR{A}$ (in that order), followed by some zero rows.

By replacing the zero rows at the bottom of \BMR{S} by some standard basis elements of $\F^{m}$, one can create an invertible matrix \BMR{S} which still satisfies \eqref{e1;5-3-24}.
\end{Cor}

\begin{proof}
Take matrices $\BMR{B}$ and $\BMR{G}$ described in theorem \ref{T2;18-1-24}. If we beef $\BMR{B}$ up to an $n\times m$ size, by attaching zero columns from the right, the resulting matrix is exactly $\RCEF{\BMR{A}}$. Append a corresponding number of zero rows to the bottom of \BMR{G} to produce \BMR{S}.

Let us settle the second claim. The rows of $\BMR{G}$ form a basis of the row space of $\BMR{A}$, and so can be extended to a basis of $\F^{m}$ through addition of some standard basis elements of $\F^{m}$.\footnote{Of course this needs to be demonstrated, but it is an easy argument.} If we replace the zero rows at the bottom of \BMR{S} with these, we produce the desired new invertible matrix $\BMR{S}$.
\end{proof}

If we apply corollary \ref{c1;19-1-24} to $\BMR{A}^{T}$, use equalities \eqref{e3;19-1-24} and \eqref{e2;19-1-24}, remembering that column space of $\BMR{A}^{T}$ and $\text{Nullspace}{\left(A\BMR{}\right)}$ are $\bullet$-complements, we obtain a corresponding claim about RREF factorization.

\begin{Cor}[RREF factorization]
Suppose that $\BMR{A}$ is an $n\times m$ non-zero matrix, and $i_{1}< i_{2}<  \cdots< i_{r}$ are the \Lime indices of the row space of $\BMR{A}$.\footnote{These \Lime indices of the row space of $\BMR{A}$ are exactly the non-\Red indices of the nullspace of $\BMR{A}$.}
Then 
\begin{equation}\label{e2;5-3-24}
\BMR{A}=\BMR{T}\circ \RREF{\BMR{A}},
\end{equation}
where the columns of $\BMR{T}$ are the columns $C_{i_{1}}, C_{i_{2}}, \ldots , C_{i_{r}}$ of $\BMR{A}$ (in that order), followed by some zero columns.

By replacing the eastern zero columns of \BMR{T} by some 
some standard basis elements of $\F^{n}$, one can create an invertible matrix \BMR{T} which still satisfies \eqref{e2;5-3-24}.
\end{Cor}

\section{Gauss-Jordan \`{a} la Gram-Schmidt}\mbox{}
The reader may have grown concerned that our approach does not expose students  to Gauss-Jordan elimination that has long been sacred in a first course on linear algebra.

We live in the contemporary world with vast computing resources and AI at our fingertips. It is counterproductive to let students think that proficiency in carrying out Gauss-Jordan elimination by hand is a valuable intellectual accomplishment.  
Proficiency in reducing a quest to a recognizable procedure, being able to use computers to carry out the calculations, and having the skills to interpret the results, is of course a different matter.

The efficiency we gain through our approach is in part due to fact that row-reduction no longer needs to play a central role in the course.
Still, Gauss-Jordan elimination has a natural place within the development of \Red\!/\Lime bases, {\emph {well before matrices are officially introduced}}. Furthermore, in this form the elimination is reminiscent of Gram-Schmidt process, or vice versa. This is something that students appreciate.

The value of including this algorithm in a first linear algebra course is in demonstrating its existence. We certainly do not suggest that students should be asked to carry out such a procedure by hand; one should use a computer for any such task.

Like Gram-Schmidt process, the procedure we present is recursive, and it relies on the repeated use of a procedure $\mathfrak{P}$ (described below) which pertains to the following question: 
when one appends an element $Y\in\F^{n}$ to the \Lime basis of a subspace \Fam{W} of $\F^{n}$, how does one generate the \Lime basis of the resulting span \Fam{V} of the new list?\footnote{Since there is nothing to do if $Y$ is already in \Fam{W}, and testing for this inclusion is easy (see \eqref{e1}\!\! ), we will focus on the case of $Y\not\in\Fam{W}$.} We consider the \Lime case because it corresponds most directly to the traditional Gauss-Jordan elimination.

Before stating a general form of the procedure, let us illustrate the ideas with a concrete and simplified example. This should convince the reader that the procedure is indeed akin to Gauss-Jordan elimination.

Let us say that some $Y\in\F^{n}$ originates with a $3$ in the $7$-th position, and we insert $Y$ into the list between $W_{4}$ and $W_{8}$:
\begin{equation*}
\resizebox{\textwidth}{!}{$\begin{array}{rcccccccccccccccccc}
W_{3}&=&(0,& 0,& \boxed{1},& 0,& \star,& \star,& \boxed{{\BM a}},& 0,& \star,& \star& \star,& \star,& 0,& \ldots,& 0,& \ldots, & \star)\\
W_{4}&=&(0,& 0,& 0,& \boxed{1},& \star,& \star,& \boxed{{\BM b}},& 0,& \star,& \star,& \star,& \star,& 0,& \ldots,& 0,& \ldots &\star)\\
Y&=&(0,& 0,& 0,& 0,&  0,& 0,& \boxed{3},& \boxed{{\BM \alpha}},& \star,& \star,& \star,& \star,& \boxed{{\BM \beta}},& \ldots,& \boxed{{\BM \zeta}},& \ldots & \star)\\
W_{8}&=&(0,& 0,& 0,& 0,& 0,& 0,& 0,& \boxed{1},& \star,& \star,& \star,& \star,& 0,& \ldots,& 0,& \ldots &\star)\\
W_{13}&=&(0,& 0,& 0,& 0,& 0,& 0,& 0,& 0,& 0,& 0,& 0,& 0,& \boxed{1},& \ldots,& 0,& \ldots & \star)\\ : \\
W_{p}&=&(0,& 0,& 0,& 0,& 0,& 0,& 0,& 0,& 0,& 0,& 0,& 0,& 0,& \ldots,& \boxed{1},& \ldots & \star)\\
\end{array}$}
\end{equation*}

Then we can  zero-out the entries of $Y$ in the $8$-th, $13$-th and $p$-th positions, and scale the new \Fam{Y} so that it originates with a $1$:
\begin{quote}
Let $\hat{Y}=\frac{1}{3}\left(Y-{\BM \alpha}W_{8}-{\BM \beta}W_{13}-\cdots-{\BM \zeta}W_{p}\right)$.  
\end{quote}
Replacing $Y$ with $\hat{Y}$ does not change the span of the list.

Then
\begin{equation*}
\resizebox{\textwidth}{!}{$\begin{array}{rcccccccccccccccccc}
W_{3}&=&(0,& 0,& \boxed{1},& 0,& \star,& \star,& \boxed{{\BM a}},& 0,& \star,& \star& \star,& \star,& 0,& \ldots,& 0,& \ldots, & \star)\\
W_{4}&=&(0,& 0,& 0,& \boxed{1},& \star,& \star,& \boxed{{\BM b}},& 0,& \star,& \star,& \star,& \star,& 0,& \ldots,& 0,& \ldots &\star)\\
\hat{Y}&=&(0,& 0,& 0,& 0,&  0,& 0,& \boxed{1},& \boxed{0},& \star,& \star,& \star,& \star,& \boxed{0},& \ldots,& \boxed{0},& \ldots & \star)\\
W_{8}&=&(0,& 0,& 0,& 0,& 0,& 0,& 0,& \boxed{1},& \star,& \star,& \star,& \star,& 0,& \ldots,& 0,& \ldots &\star)\\
W_{13}&=&(0,& 0,& 0,& 0,& 0,& 0,& 0,& 0,& 0,& 0,& 0,& 0,& \boxed{1},& \ldots,& 0,& \ldots & \star)\\ : \\
W_{p}&=&(0,& 0,& 0,& 0,& 0,& 0,& 0,& 0,& 0,& 0,& 0,& 0,& 0,& \ldots,& \boxed{1},& \ldots & \star)\\
\end{array}$}
\end{equation*}

Next we replace $W_{3}$ and $W_{4}$ by $\overline{W_{3}}=W_{3}-a\hat{Y}$ and $\overline{W_{4}}=W_{4}-b\hat{Y}$ respectively, to produce the list
\begin{equation*}
\resizebox{\textwidth}{!}{$\begin{array}{rcccccccccccccccccc}
\overline{W_{3}}&=&(0,& 0,& \boxed{1},& 0,& \star,& \star,& \boxed{0},& 0,& \star,& \star& \star,& \star,& 0,& \ldots,& 0,& \ldots, & \star)\\
\overline{W_{4}}&=&(0,& 0,& 0,& \boxed{1},& \star,& \star,& \boxed{0},& 0,& \star,& \star,& \star,& \star,& 0,& \ldots,& 0,& \ldots &\star)\\
\hat{Y}&=&(0,& 0,& 0,& 0,&  0,& 0,& \boxed{1},& \boxed{0},& \star,& \star,& \star,& \star,& \boxed{0},& \ldots,& \boxed{0},& \ldots & \star)\\
W_{8}&=&(0,& 0,& 0,& 0,& 0,& 0,& 0,& \boxed{1},& \star,& \star,& \star,& \star,& 0,& \ldots,& 0,& \ldots &\star)\\
W_{13}&=&(0,& 0,& 0,& 0,& 0,& 0,& 0,& 0,& 0,& 0,& 0,& 0,& \boxed{1},& \ldots,& 0,& \ldots & \star)\\ : \\
W_{p}&=&(0,& 0,& 0,& 0,& 0,& 0,& 0,& 0,& 0,& 0,& 0,& 0,& 0,& \ldots,& \boxed{1},& \ldots & \star)\\
\end{array}$}
\end{equation*}
We have not altered the span, and the new list is the \Lime basis of the span of $W_{3}, W_{4}, W_{8}, W_{13}, \ldots, W_{p}, Y$.
\footnote{
Of course, we could have done the same type of thing with right-standard bases.}
\BSK

Here is the general {``Procedure $\mathfrak{P}$''} for the case when the \Lime signature of a subspace \Fam{W} of $\F^{n}$ is $\left(3,4,8,13,\ldots,p\right)$, and 
$W_{3}, W_{4}, W_{8}, W_{13}, \ldots, W_{p}$ are the corresponding \Lime basis elements of \Fam{W}. 

Start by replacing $Y$ with $\hat{Y}=Y-Y(3)W_{3}-Y(4)W_{4}-Y(8)W_{8}-\cdots-Y(p)W_{p}$. The span of $W_{3}, W_{4}, W_{8}, W_{13}, \ldots, W_{p}, \hat{Y}$ is still \Fam{V}, i.e. the span of $W_{3}, W_{4}, W_{8}, W_{13}, \ldots, W_{p}, Y$. 

Clearly $\hat{Y}$ has zero entries in every \Lime position of \Fam{W}. In particular $Y\in \Fam{W}$ if and only if $\hat{Y}$ is the zero vector. 

We only need to deal with the case when $\hat{Y}$ is not the zero vector, and in that case  
the dimension of \Fam{V} is one larger than that of \Fam{W}. Hence \Fam{V} gains a new \Lime index in addition to those inherited from \Fam{W}: the originating index $i_{o}$ of $\hat{Y}$. After scaling, we may assume that $\hat{Y}$ originates with a $1$ (and still has zero entries in all \Lime positions for \Fam{W}).
This newly minted $\hat{Y}$ is the \Lime basis element of \Fam{V} corresponding to the index $i_{o}$.

The original \Lime basis elements of \Fam{W} corresponding to the \Lime indices that exceed $i_{o}$ also serve as \Lime basis elements of \Fam{V}.

The \Lime basis elements of \Fam{W} corresponding to the indices preceding $i_{o}$ would also be \Lime basis elements of \Fam{V}, were it not for their potentially non-zero $i_{o}$-th entries. So we subtract from them appropriate multiples of $\hat{Y}$ to zero out these entries, and, voil\`a,  we arrive at the \Lime basis of \Fam{V} (Observation \ref{Obs1;28-5-25}).

\BSK
By using the procedure $\mathfrak{P}$ recursively, we can produce the \Lime basis $Y_{1}, Y_{2},  \ldots, Y_{p}$  of the span of a list $O\neq X_{1}, X_{2},  \ldots, X_{k}$ in $\F^{n}$, with well-aligned intermediate spans, \`{a} la Gram-Schmidt process. 

Start by scaling $X_{1}$ so that it originates with a $1$. This is $Y_{1}$, and the case $k=1$ is done. Then use the procedure $\mathfrak{P}$ recursively on $k$, by appending $X_{j+1}$ to the already constructed \Lime basis $Y_{1}, Y_{2},  \ldots, Y_{l}$ of the span of $X_{1}, X_{2},  \ldots, X_{j}$.

Of course this algorithm can be used to construct the \Lime basis of the row space of $\BMR{A}$, i.e. to construct $\RREF{\BMR{A}}$.

\section{Sub-Terminal Entries, Subspace Truncations And Possible Configurations}\mbox{}

This section is motivated by a  natural question (which turns out to have a pretty answer): ``what configurations of {\color{red} red} and {\color{lime!50!black} lime} indices are possible?''

We will work our way towards the answer by asking a simpler question:
``what can be said about the positions of the second-to-last non-zero entry in the elements of a subspace \Fam{W} of $\F^{n}$?'' 

Let us say that $Z\in\Fam{W}$ \textit{\textbf{sub-terminates}} in the $j$-th position if its second-to-last non-zero entry is in that position. If an element has no second-to-last non-zero entry, its sub-terminal index is $0$.

When $X$ and $Y$ are distinct elements of \Fam{W} that terminate in the $i$-th position, after scaling we can assume that their terminating non-zero entries (in the $i$-th position) are both $1$. Hence $X-Y$ is an element of \Fam{W} that terminates in some {\color{red} red} $k$-th position, with $k<i$. In particular, $X$ and $Y$ coincide past the $k$-th position, but not in the $k$-th position. Consequently, the sub-terminal index of at least one of $X, Y$ is at least $k$ (and is less than $i$). If it is strictly between $k$ and $i$, then it is the sub-terminal index for both $X$ and $Y$. This shows that when the sub-terminal entries of $X$ and $Y$ occur in non-{\color{red} red} positions, they occur in the same position.

This shows that  
\MD{there is at most one non-{\color{red} red} position in which the elements of \Fam{W} that terminate in the $i$-th position can sub-terminate.}
The sub-terminal index $j_{_{i}}$ of a {\color{red} red}-basic element $W_{i}$ of \Fam{W} \eqref{T3} is not {\color{red} red}, and so it is the only non-{\color{red} red} sub-terminal index for elements of \Fam{W} that terminate in the $i$-th position. Hence we can identify all of the sub-terminal indices of the elements of \Fam{W} that terminate in the $i$-th position. \MD{These indices are: $j_{_{i}}$,  and all {\color{red} red} indices between $j_{_{i}}$ and $i$. \footnote{To reach this conclusion one makes use of the {\color{red} red}-basic elements for indices preceding $i$.}\BSK}

It will be handy to designate the $n$ indices/positions of a subspace \Fam{W} of $\F^{n}$ by letters $\rho, \lambda,  \beta$ and $\nu$ as follows:
\begin{equation*}
\left\{\begin{array}{rl}
\rho:& \text{{\color{red}{\bf r}ed}, but not {\color{lime!50!black} lime}}\\
\lambda:& \text{{\color{lime!50!black}{\bf l}ime}, but not {\color{red} red}}\\
\beta:& \text{{\bf b}oth {\color{red} red} and {\color{lime!50!black} lime}}\\
\nu:& \text{{\bf n}either {\color{red} red} nor {\color{lime!50!black} lime}}
\end{array}\right.
\end{equation*}\MD{The last position for a subspace \Fam{W} of $\F^{n}$ is a $\beta$-position exactly when \Fam{W} contains the standard basis element $E_{n}$. 

The last position is a $\nu$-position if and only if all elements of \Fam{W} have zero as their last entry. 

The last position of a subspace can never be a $\lambda$-position, and the first position - can never be a $\rho$-position. }
We can transform $\F^{n}$, with $n>1$, into $\F^{n-1}$, by simply dropping off the last entry of the $n$-tuples to produce $(n-1)$-tuples. In the process subspaces of $\F^{n}$ get transformed into subspaces of $\F^{n-1}$. For simplicity, let us refer to this as a ``truncation of $\F^{n}$ (and its subspaces) from the right.''

\begin{Thm}[Effect of truncation on {\color{red} red}/{\color{lime!60!black} lime} indices]\label{T1}
The act of a truncation of \Fam{W} from the right to produce a subspace \Fam{U} has the following effect on the special positions:
\begin{itemize}
\item If the last position of \Fam{W} is a $\beta$-position or a $\nu$-position, then the indices of \Fam{U} have exactly the same designations as they had for \Fam{W}.\SSK
\I If the last position of \Fam{W} is a  $\rho$-position, then one of the indices for \Fam{U} gains a {\color{red} red} status (that it did not have for \Fam{W}), and status of all other indices remain unchanged. So, either some $\nu$-index of \Fam{W} becomes a $\rho$-index for \Fam{U}, or a $\lambda$-index for \Fam{W} becomes a $\beta$-index for \Fam{U}. 
\end{itemize}
\end{Thm}

\begin{proof}The only possible change in status is that sub-terminal positions of the elements of \Fam{W} that terminate in the last ($n$-th) position, become terminal positions for the corresponding elements of \Fam{U}.

As we have observed, such sub-terminal positions are either the sub-terminal position of the {\color{red} red}-basic element of \Fam{W} corresponding to the index $n$ (if $n$ happens to be a {\color{red} red} index) or are the positions corresponding to the {\color{red} red} indices that are less than $n$. \qedhere
\end{proof}

Now we are ready to answer our original question. Suppose that $n$ positions are marked by (perhaps repeated) symbols $\rho, \lambda, \beta, \nu$, in some fashion. Is there a subspace \Fam{W} of $\F^{n}$ for which this is exactly its ``signature'', in the sense described above?

The answer is certainly ``not always''. As we have already noted, the last position cannot be a $\lambda$-position. So, which signatures are possible?

We begin with an example that provides a partial answer to the question by demonstrating a sufficient condition.

Let
\begin{equation*}
\Fam{W}=\Set{\left(0, a_{1}, a_{2}, a_{3}, a_{4}, 0, a_{1}, a_{5}, a_{6}, 0, a_{4}, a_{6},  a_{7},  a_{7}, a_{8}, a_{9}, a_{8}, a_{10}  \right)}{a_{i}\in\F}.
\end{equation*}

Then \Fam{W} is a subspace of $\F^{18}$ with the following signature:
\begin{equation*}
\left(\nu, \lambda, \beta, \beta, \lambda, \nu, \rho, \beta, \lambda, \nu, \rho, \rho, \lambda, \rho, \lambda, \beta, \rho,  \beta \right).
\end{equation*}

From this example it is quite obvious that the locations of the $\nu$-positions and $\beta$-positions are generally unrestricted. Furthermore, if the locations of the $\lambda$-positions and $\rho$-positions correspond to the positions of ``matched'' left/right parentheses, the whole configuration is achievable as a signature of a subspace of $\F^{n}$ of the type shown. This describes a sufficient condition, but is it necessary? We claim that it is.

Any configuration of matched left/right parentheses has the property that to the left of any left parenthesis there are more right parentheses than there are left ones. In fact, apart from stating that there are as many left parentheses as there are right parentheses, this property identifies all possible configurations of matched left/right parentheses. 

\begin{Thm}[All Possible Configurations]
A $\rho/\lambda/\beta/\nu$-configuration of $n$ positions represents a signature of a subspace of $\F^{n}$ if and only if the number of $\lambda$-positions equals the number of $\rho$-positions, and to the right of every $\lambda$-position there are more $\rho$-positions than $\lambda$-positions.
\end{Thm}

\begin{proof}
We focus on the forward implication, since the validity of the reverse implication has already been discussed. Our proof is inductive on $n$. The claim is trivially true when $n=1$. If it were not universally true, there would be a smallest $n_{o}>1$ for which the claim fails. So, we have some subspace $\Fam{W}_{o}$ of $\F^{n_{o}}$ whose signature does not have the stated property. The signature of the right truncation $\Fam{U}_{o}$ of $\Fam{W}_{o}$ would enjoy the described property by the minimality of $n_{o}$.

This means (by theorem \ref{T1}) that the last position of $\Fam{W}_{o}$ would have to be a $\rho$-position. Yet in that case the signature of $\Fam{W}_{o}$ could be recovered from the signature of $\Fam{U}_{o}$ by turning a $\rho$-position into a $\nu$-positions or a $\beta$-position into a $\lambda$-position, and attaching a new  $\rho$-position to the end. Such a transformation could not produce a $\lambda$-position that has no more $\rho$-positions to the right of it than $\lambda$-positions, in contradiction to hypothesis about $\Fam{W}_{o}$.
\end{proof}
It is not hard to see that a feasible $\rho/\lambda/\beta/\nu$-signature does not always describe a unique subspace. How different can two subspaces be if they have identical  $\rho/\lambda/\beta/\nu$-signatures? Note that the persistent zero position may not be the same. Subspaces 
\begin{equation*}
\Fam{Z}=\Set{\left(0, a_{1}, a_{2}, a_{3}, a_{4}, a_{4}, a_{1}, a_{5}, a_{6}, a_{6}, a_{4}, a_{6},  a_{7},  a_{7}, a_{8}, a_{9}, a_{8}, a_{10}  \right)}{a_{i}\in\F}
\end{equation*}
and
\begin{equation*}
\Fam{X}=\Set{\left(0, a_{1}, a_{2}, a_{3}, a_{4}, a_{1}, a_{1}, a_{5}, a_{6}, a_{4}, a_{4}, a_{6},  a_{7},  a_{7}, a_{8}, a_{9}, a_{8}, a_{10}  \right)}{a_{i}\in\F}
\end{equation*}
of $\F^{18}$ both have the same signature as the subspace \Fam{W} we had presented earlier. 

Of course if a signature begins with a $\nu$, this forces the first position of a subspace to be a zero position. Analogous claim is true when $\nu$ is the last symbol of a signature.

\subsection{One More Thing}
So, before we wrap things up, there is one more question that seems too natural to be ignored, and has a simple answer.

\begin{quote}
If a subspace \Fam{Z} of $\F^{n}$, restricted to some $k$ positions, presents as $\F^{k}$,  can one apply a permutation matrix to \Fam{Z} and have those $k$ positions become \Red positions of the resulting subspace? 
\end{quote}

The answer is affirmative. Take a permutation matrix that maps those $k$ positions of $\F^{n}$ to the last $k$ positions (with indices $\scriptstyle n-(k-1),\ \ldots,\ n-1,\ n$), and apply it to \Fam{Z} to produce a subspace $\Fam{Z}_{_{o}}$. Since, restricted to these last $k$ positions, $\Fam{Z}_{_{o}}$ presents as $\F^{k}$, each of these positions is indeed a \Red position for $\Fam{Z}_{_{o}}$. 

If $k$ happens to be the dimension of \Fam{Z}, then the last $k$ positions of $\Fam{Z}_{_{o}}$ are all of its \Red positions. 
In this case, the elements of $\Fam{Z}$ presenting as the standard basis of $\F^{k}$ on the original $k$ positions, are uniquely determined, and are sent to the \Red basis of $\Fam{Z}_{_{o}}$ by the permutation matrix.

\section{Concluding Remarks}

According to John Baez' insight \cite{Baez}, a good mathematical article should not end abruptly, and we concur. 

Originally we were motivated by a search for a palatable argument showing the uniqueness of RREF for the students in the first linear algebra course at our home institution. At the same time we had grown weary of  our students fixating on row reduction and the RREF, making these a knee-jerk response to almost any linear-algebraic question. Given the  centrality of the orthogonality-related notions in contemporary applications, this just won't do.

Furthermore, it is counterproductive to let students think that proficiency in carrying out Gauss-Jordan elimination or a determinant calculation by hand is a sign of valuable intellectual accomplishment. Present day professionals perform such procedures by hand just about as often as they extract square roots by hand or graph complicated functions with pencils on paper. Proficiency in reducing a quest to Gauss-Jordan elimination (or a determinant calculation), along with an ability to use computing systems, and having the skills to interpret the results, is of course a different matter.

It has to be mentioned that we have only tested ``in a jiffy'' approach at our home institution, which benefits from having strong and motivated students. It works well and is indeed efficient. The approach is most effective when the emphasis is placed  on comprehension of the ideass, and the general proofs are scaled down to more visual and concrete arguments that can be played with in real time.

It is our hope that the \Red/\Lime bases tack can be adopted more broadly, and that it finds a regular place in the toolbox of modern mathematical educators.\Vm{3}

\textit{Acknowledgements: We are very grateful to Professors Fernando Gouvea and Scott Taylor, who had kindly read the first version of the article. Their insightful suggestions prompted a significant improvement, as did the blog post by John Baez which Fernando Gouvea had brought to my attention. We would also be remiss not to thank those Colby students whom we had subjected to early versions of the ideas in this article. Their questions and comments contributed greatly to a refinement of the approach.}

\bibliographystyle{vancouver}
\bibliography{MyBibl}

\mbox{}
\vspace{3em}
\vfill

\end{document}